\documentclass[12pt]{article}

\usepackage{amsfonts}

\def\be{\begin{equation}}
\def\ee{\end{equation}}
\def\bea{\begin{eqnarray}}
\def\eea{\end{eqnarray}}

\def\1{\'{\i}}                         
\def\R{{\mathbb R}}

\def\a#1{a_{#1}}
\def\b#1{b_{#1}}
\def\c#1{c_{#1}}

\def\ap{A_+}
\def\am{A_-}
\def\bp{B_+}
\def\bm{B_-}
\def\aa{N}
\def\bb{M}

\def\aap{{\hat a}_+}
\def\aam{{\hat a}_-}

\begin{document}
\

\begin{center}

{\LARGE{\bf{Two-photon algebra deformations}}}

\vskip1cm

{\sc Preeti Parashar, Angel Ballesteros and Francisco J. Herranz}
\vskip0.5cm

{\it Departamento de F\1sica, Universidad de Burgos}

{\it 09001 Burgos, Spain}

\end{center}
\vskip1cm

\begin{abstract}
\noindent
In order to obtain a classification of all possible  quantum
deformations of the two-photon algebra $h_6$, we introduce its
corresponding  general Lie bialgebra, which is a  coboundary one.  Two
non-standard quantum deformations  of $h_6$, together with  their
associated quantum  universal $R$-matrix, are presented; each of them
contains either a quantum   harmonic oscillator subalgebra  or a quantum
$gl(2)$ subalgebra.  One-boson representations for these quantum
two-photon algebras are derived and translated into Fock--Bargmann
realizations. In this way, a systematic  study of `deformed' states of
light in quantum optics can be developed.
\end{abstract}

\vskip1cm

\section{Introduction}

The two-photon Lie algebra $h_6$  is generated by the
operators $\{\aa,\ap,\am,\bp,\bm,\break \bb\}$ endowed with the following
commutation rules \cite{Gil}:
\be
\begin{array}{lll}
 [\aa,\ap]=\ap \quad  &[\aa,\am]=-\am \quad &[\am,\ap]=\bb \cr
[\aa,\bp]=2\bp \quad  &[\aa,\bm]=-2\bm \quad &[\bm,\bp]=4\aa+2\bb
\cr
 [\ap,\bm]=-2\am \quad  &[\ap,\bp]=0 \quad
&[\bb,\,\cdot\,]=0\cr
 [\am,\bp]=2\ap \quad  &[\am,\bm]=0 .\quad
& \cr
\end{array}
\label{aa}
\ee
The  Lie algebra $h_6$ contains several remarkable Lie subalgebras: the
Heis\-enberg--Weyl algebra
$h_3$ spanned by $\{\aa,\ap,\am\}$, the harmonic
 oscillator algebra $h_4$ with
generators 
$\{\aa,\ap,\am,\bb\}$,  and the $gl(2)$ algebra
generated by $\{\aa,\bp,\bm,\bb\}$. Note that  $gl(2)$ is isomorphic to 
a trivially extended ${sl}(2,\R)$ algebra (the central  extension
$\bb$   can be absorbed by redefining $\aa\to\aa+\bb/2$). Hence we have
the following embeddings
\be
h_3\subset h_4\subset h_6\qquad sl(2,\R)\subset gl(2)\subset h_6 .
\label{ab}
\ee

Representations of the two-photon algebra can be used to  generate a
large zoo of squeezed  and coherent states for (single mode) one- and
two-photon processes which have been  analysed in \cite{Gil,Brif}. 
In particular, if the generators $\aam$, $\aap$ close a boson algebra
\be
[\aam,\aap]=1 ,
\label{ac}
\ee
then a one-boson representation of $h_6$ reads
\be
\begin{array}{lll}
\aa=\aap\aam&\qquad \ap=\aap&\qquad \am=\aam\cr
\bb=1&\qquad \bp=\aap^2&\qquad \bm=\aam^2  .
\end{array}
\label{ad}
\ee
This realization shows that one-photon processes  are algebraically
encoded within the subalgebra $h_4$, while $gl(2)$ contains the
information concerning two-photon dynamics.
 When the operators
$\aam$, $\aap$ act in the usual way on the number states Hilbert space
spanned by $\{|m\rangle\}_{m=0}^{\infty}$, {\it i.e.},
\be
\aam |m\rangle = \sqrt{m}\,|m-1\rangle   \qquad 
\aap |m\rangle = \sqrt{m+1}\,|m+1\rangle ,
\label{ae}
\ee
the action of   $h_6$  on these states becomes
\be 
\begin{array}{ll}
N|m\rangle =m |m\rangle &\quad M|m\rangle=|m\rangle\\
 A_+|m\rangle= \sqrt{m+1}\, |m+1\rangle  
 &\quad B_+|m\rangle= \sqrt{(m+1)(m+2)}\,|m+2\rangle\\
A_-|m\rangle = \sqrt{m}\,|m-1\rangle  &\quad 
B_-|m\rangle=  \sqrt{m(m-1)}\,|m-2\rangle.
\end{array}
\label{af}
\ee

The one-boson realization (\ref{ad}) can be translated  into a
Fock--Bargmann representation \cite{FB} by  setting $\aap\equiv \alpha$
and $\aam\equiv
\frac d{d\alpha}$. Thus  the $h_6$ generators act in the
Hilbert space of entire   analytic functions
$f(\alpha)$ as linear differential operators:
\bea
&&\aa=\alpha \frac{d}{d\alpha} \qquad \ap=\alpha \qquad 
\am=\frac{d}{d\alpha}\nonumber\\ 
&&\bb=1 \qquad\quad \bp=\alpha^2 \qquad
\bm=\frac{d^2}{d\alpha^2} .
\label{ag}
\eea
The two-photon algebra eigenstates \cite{Brif} are given by the
analytic eigenfunctions that fulfil
\be
(\beta_1\aa + \beta_2 \bm + \beta_3\bp+\beta_4\am+\beta_5\ap)
f(\alpha)=\lambda f(\alpha).
\label{ah}
\ee
In the Fock--Bargmann
representation (\ref{ag}), the
following differential equation is deduced from (\ref{ah}):
\be
\beta_2 \frac{d^2f}{d\alpha^2}+(\beta_1\alpha  +\beta_4)\frac{d
f}{d\alpha} +(\beta_3\alpha^2+\beta_5\alpha -\lambda)f=0 
\label{ai}
\ee
where $\beta_i$ are arbitrary complex coefficients and  $\lambda$ is a
complex eigenvalue. The solutions of this equation  (provided a
suitable normalization is imposed) give rise to the two-photon
coherent/squeezed states \cite{Brif}. One- and two-photon coherent and
squeezed states corresponding to the subalgebras $h_4$ and
$gl(2)$ can be derived from equation (\ref{ai}) by
setting $\beta_2=\beta_3=0$ and  $\beta_4=\beta_5=0$, respectively.

The  prominent role that the two-photon Lie algebra  plays in relation
with squee\-zed and coherent states motivates the extension of the Lie
bialgebra classifications already perfor\-med for its 
subalgebras $h_3$ \cite{heis},  $h_4$ \cite{osc} and $gl(2)$
\cite{gl}, since  the  two-photon bialgebras    would constitute the
underlying structures of any   further  quantum deformation  whose
representations  could be physically interesting in the field of quantum
optics. Thus in the next section we present such a classification for
the two-photon bialgebras. The remaining sections are devoted to show
how quantum two-photon deformations provide  a starting point    in the
analysis of `deformed' states of light.


\section{The two-photon Lie bialgebras}

The essential point in this contribution is the fact that any quantum
deformation of a given Lie algebra can be characterized (and sometimes
obtained) through the associated Lie bialgebra.

Let us first recall that a {\it Lie bialgebra} $(g,\delta)$ is a Lie
algebra $g$ endowed with a linear  map $\delta:g\to g\otimes g$ called 
the {\it cocommutator} that fulfils two conditions \cite{CP}:

\noindent i) $\delta$ is a 1-cocycle, {\it i.e.},
\be
\delta([X,Y])=[\delta(X),\, 1\otimes Y+ Y\otimes 1] +
[1\otimes X+ X\otimes 1,\, \delta(Y)]  \quad \forall X,Y\in
g.
\label{ba}
\ee
\noindent ii) The dual map $\delta^\ast:g^\ast\otimes g^\ast \to
g^\ast$ is a Lie bracket on $g^\ast$.

A Lie bialgebra $(g,\delta)$ is called a  {\it coboundary} Lie
bialgebra  if there exists an element $r\in g\wedge g$ called the {\it
classical
$r$-matrix}  such that
\be
\delta(X)=[1\otimes X + X \otimes 1,\,  r]  \qquad  \forall X\in
g.
\label{bb}
\ee
Otherwise the Lie bialgebra is a {\it non-coboundary} one.

There are two types of coboundary Lie bialgebras $(g,\delta (r))$:

\noindent
i)  {\em Non-standard} (or triangular): The $r$-matrix is a
skewsymmetric solution of the classical Yang--Baxter equation (YBE):
\be
[[r,r]]=0 ,
\label{bc}
\ee
where $[[r,r]]$ is the Schouten bracket defined by
\be
[[r,r]]:=[r_{12},r_{13}] + [r_{12},r_{23}] + 
[r_{13},r_{23}].
\label{bd}
\ee
If $r=r^{i j} X_i\otimes X_j$, we have denoted
$r_{12}=r^{i j} X_i\otimes X_j\otimes 1$,
$r_{13}=r^{i j} X_i\otimes 1\otimes X_j$ and
$r_{23}=r^{i j} 1\otimes X_i\otimes  X_j$.

\noindent
ii) {\em Standard} (or quasitriangular): The $r$-matrix is a
skewsymmetric solution of the\break modified classical  YBE:
\be
[X\otimes 1\otimes 1 + 1\otimes X\otimes 1 +
1\otimes 1\otimes X,[[r,r]]\, ]=0  \quad \forall X\in g.
\label{be}
\ee

\subsection{The general solution}

Now we proceed to introduce {\it all}  the Lie bialgebras associated
to $h_6$. Recently a classification of    all Schr\"odinger bialgebras
have been obtained \cite{sch} showing that all of them have a coboundary
character. Therefore we can make use of the isomorphism between the 
  Schr\"odinger and two-photon algebras  \cite{tufotona} in order to
`translate' the results of the former in terms of the latter.

The most general two-photon classical $r$-matrix,  $r\in h_6\wedge h_6$,  
 depends on 15 (real) coefficients:
\bea
&& r=\a1 N\wedge A_+ + \a2 N\wedge B_+  +
\a3 A_+\wedge M  \cr
&&\qquad + \a4 B_+\wedge M +
\a5 A_+\wedge B_+ + \a6 A_+\wedge B_-\cr
&&\qquad + \b1 N\wedge A_-  + \b2 N\wedge B_-  +
\b3 A_-\wedge M  \cr
&&\qquad + \b4 B_-\wedge M +
\b5 A_-\wedge B_- + \b6 A_-\wedge B_+\cr
&&\qquad +\c1 N\wedge M + \c2 A_+\wedge A_- +
\c3  B_+\wedge  B_-  
\label{bf}
\eea 
which are subjected to 19 equations that we group into three sets:
\be 
\begin{array}{l}
2 \a6^2 - \a6 \b1 + 3 \a1 \b5 + 2 \b5 \b6=0\cr
\a2 \a3-2 \a1 \a4+2 \a4 \b6-3 \a5 \c1-\a5 \c2-2 \a5 \c3=0\cr
\a1 \a2-2 \a2 \b6-4 \a5 \c3=0\cr
\a5 \b1-\a1 \b6+2 \a2 \c1+2 \a2 \c3+4 \a4 \c3=0\cr
2 \a2 \a6+4 \a4 \a6-2 \a4 \b1-2 \a5 \b2+2 \a2 \b3-
4 \a5 \b4+\a1 \c1+\a1 \c2=0\cr
3 \a1 \b2+2 \a2 \b5+4 \a6 \c3-2 \b1 \c3=0\cr
\a3 \b2+2 \a1 \b4+2 \a4 \b5+\a6 \c1-\a6 \c2-2 \a6 \c3-2 \b3 \c3=0\cr
3 \a2 \b5+\b2 \b6+2 \a6 \c3=0
\end{array}
\label{bbf}
\ee  
\be 
\begin{array}{l}
2 \b6^2 - \b6 \a1  + 3 \b1 \a5   + 2 \a5 \a6 =0   \cr
\b2 \b3 -2 \b1 \b4 +2  \b4 \a6 -3 \b5 \c1+\b5 \c2+2 \b5 \c3=0\cr
\b1 \b2 - 2 \b2 \a6 +4 \b5 \c3=0\cr
\b5 \a1  - \b1 \a6 +2 \b2 \c1-2 \b2 \c3-4 \b4 \c3=0\cr
2 \b2 \b6 +4 \b4 \b6-2 \b4 \a1 -2 \b5 \a2 + 2 \b2 \a3 -
4 \b5 \a4 +\b1 \c1-\b1 \c2=0\cr
3 \b1 \a2 +2 \b2 \a5-4 \b6 \c3  +2 \a1 \c3=0\cr
\b3 \a2 + 2 \b1 \a4 +2 \b4 \a5 +\b6 \c1+\b6 \c2+2 \b6 \c3+2 \a3 \c3=0\cr
3\b2 \a5 + \a2 \a6-2 \b6 \c3=0 
\end{array}
\label{bbg}
\ee  
\be 
\begin{array}{l}
\a2 \b2 + \c3^2=0\cr
2 \a2 \b4 + 2  \a4 \b2  - \a5 \b5 + \a6 \b6 - 2 \c3^2=0\cr
\a1 \b1 + \a1 \a6+ \b1 \b6 +2 \a5 \b5-2 \a6 \b6=0 .
\end{array}
\label{bbh}
\ee

The classical $r$-matrix (\ref{bf}) satisfies the modified
classical YBE and its  Schouten bracket  reads
\be
[[r,r]]= ( \a1 \b3 + \a3 \b1 + 2 \a3 \a6  + 2 \b3 \b6
  - 2 \a5 \b5 + 2 \a6 \b6  - \c2^2)\ap\wedge\am\wedge\bb .
\label{bi}
\ee
Hence we obtain an additional equation which allows  us to distinguish
between non-standard and standard classical $r$-matrices:
\be 
\begin{array}{ll}
\mbox{Non-standard:}&\quad \a1 \b3 + \a3 \b1 + 2 \a3 \a6  + 2 \b3 \b6
  - 2 \a5 \b5 + 2 \a6 \b6  - \c2^2=0\cr
\mbox{Standard:}&\quad \a1 \b3 + \a3 \b1 + 2 \a3 \a6  + 2 \b3 \b6
  - 2 \a5 \b5 + 2 \a6 \b6  - \c2^2\ne 0 .
\end{array}
\label{bj}
\ee

On the other hand, the following automorphism of $h_6$
\be 
\begin{array}{lll}
N\to -N&\qquad A_+\to -A_- &\qquad A_-\to -A_+\cr
M\to -M&\qquad B_+\to -B_-&\qquad B_-\to -B_+
\end{array}
\label{bg}
\ee 
interchanges the roles
of  $A_+$ with $A_-$, and $B_+$ with $B_-$. This map can be
also implemented at a bialgebra level by introducing a suitable
transformation of the  parameters $a$'s, $b$'s and $c$'s given by
\bea
&&a_i\to b_i\qquad\ b_i\to a_i\qquad i=1,\dots,6\cr
&&\c1\to \c1\qquad \c2\to -\c2\qquad\c3\to -\c3 .
\label{bh}
\eea 
Notice that the maps (\ref{bg}) and (\ref{bh}) leave   the general
classical
$r$-matrix (\ref{bf}), the equations (\ref{bbh}) and  the Schouten
bracket (\ref{bi}) invariant, while they interchange the sets of
equations (\ref{bbf}) and (\ref{bbg}).

As all the two-photon Lie bialgebras are coboundary ones, that is, they
come from  the classical $r$-matrix (\ref{bf}), the corresponding general
cocommutator can be derived from (\ref{bb}):
\be 
\begin{array}{l}
 \delta(\aa)=\a1 \aa\wedge\ap + 2\a2 \aa\wedge \bp +
 \a3 \ap\wedge \bb + 2 \a4
\bp\wedge\bb + 3 \a5 \ap \wedge \bp \cr
 \qquad\quad -\b1 \aa\wedge\am - 2\b2 \aa\wedge \bm -
 \b3 \am\wedge \bb - 2 \b4
\bm\wedge\bb - 3 \b5 \am \wedge \bm \cr
\qquad\quad - \a6 \ap\wedge \bm + \b6 \am\wedge \bp \nonumber\\[6pt]
 \delta(\ap)=(2\a6 + \b1) \am\wedge \ap + \a2 \bp\wedge \ap 
+ \b2(\bm\wedge \ap
+2 \am\wedge \aa)\cr
\qquad\quad -\b1 \aa\wedge \bb - 2 \b4 \am\wedge\bb +  \b5 \bm\wedge\bb
+ \b6 
\bp \wedge \bb\cr
\qquad\quad - (\c1 + \c2)\ap\wedge \bb + 2\c3 \am\wedge\bp \\[6pt]
 \delta(\am)=-(2\b6 + \a1) \ap\wedge \am - \b2 \bm\wedge \am  -
\a2(\bp\wedge \am +2  \ap\wedge \aa)\cr
\qquad\quad +\a1 \aa\wedge \bb + 2 \a4 \ap\wedge\bb -  \a5 \bp\wedge\bb
- \a6 
\bm \wedge \bb\cr
\qquad\quad + (\c1 - \c2)\am\wedge \bb + 2\c3 \ap\wedge\bm  \\[6pt]
\delta(\bp)= 4 \c3 \aa\wedge\bp   + 2 (\a1 - \b6)\ap\wedge \bp + 2 \b1
\am\wedge \bp + 2 \b2 \bm\wedge \bp\cr
\qquad\quad  +2(2\a6 - \b1) \aa\wedge\ap + 2 \b5 (2\aa \wedge \am
-   \ap \wedge \bm -   \am \wedge \bb)\cr
\qquad\quad
-2(\b2 + 2\b4)\aa\wedge \bb - 2(\a6+\b3)\ap \wedge \bb 
-2(\c1+\c3)\bp\wedge\bb
\\[6pt]
\delta(\bm)= 4 \c3 \aa\wedge\bm  -  2 (\b1 - \a6)\am\wedge \bm - 2 \a1
\ap\wedge \bm - 2 \a2 \bp\wedge \bm\cr
\qquad\quad  -2(2\b6 - \a1) \aa\wedge\am - 2 \a5 (2\aa \wedge \ap
- \am \wedge \bp - \ap \wedge \bb)\cr
\qquad\quad
+2(\a2 + 2\a4)\aa\wedge \bb + 2(\b6+\a3)\am \wedge \bb 
+2(\c1-\c3)\bm\wedge\bb\\[6pt]
\delta(\bb)= 0 .
\end{array}
\label{bk}
\ee
The bialgebra automorphism  defined by the maps  (\ref{bg})
and (\ref{bh}) also interchanges the cocommutators 
$\delta(A_+)\leftrightarrow\delta(A_-)$ and
$\delta(B_+)\leftrightarrow\delta(B_-)$ leaving  $\delta(N)$ and
$\delta(M)$ invariant.

\subsection{The two-photon Lie bialgebras with two primitive
generators}

We have just shown that  in all two-photon bialgebras  the
central generator $M$ is always {\em primitive}, that is, its
cocommutator vanishes. In general,  primitive generators
determine   the physical properties of the corresponding
quantum deformations. 
Therefore we study  now those particular two-photon bialgebras
with {\it one} additional  primitive generator $X$ (besides
$M$). Furthermore, the restrictions implied by the condition 
$\delta(X)=0$ rather simplify the equations
(\ref{bbf})--(\ref{bbh}). 

Due to the equivalence
 $+\leftrightarrow -$ defined by the maps (\ref{bg}) and
(\ref{bh})  it suffices to restrict our study to  three
types of bialgebras: those with either
$N$, $A_+$ or $B_+$ primitive.

\noindent
$\bullet$ {\it Type I: $N$ primitive}.
The condition $\delta(N)=0$  leaves three free parameters $\c1$,
$\c2$ and $\c3$,  all others being equal to zero. 
The equations (\ref{bbf})--(\ref{bbh}) imply that $\c3=0$. 
The Schouten bracket  reduces to $[[r,r]]=-\c2^2\,  A_+\wedge
A_- \wedge M$, then we have a  standard subfamily with  two-parameters
$\{\c1,\c2\ne 0\}$ and a  non-standard subfamily with $\c1$ as the  only
free parameter; they read 
\bea 
&&\mbox{Standard subfamily:}\quad  \c1,\ \c2\ne 0 \nonumber\\[2pt]
&&r=\c1 N\wedge M+\c2   A_+\wedge A_-  \nonumber\\[2pt]
&&\delta(N)=0\qquad \delta(M)=0\label{ca} \\
&&\delta(A_+)=-(\c1+\c2) A_+\wedge M\qquad 
\delta(A_-)=(\c1-\c2) A_-\wedge M\nonumber\\[2pt]
&&\delta(B_+)=-2\c1  B_+\wedge M\qquad 
\delta(B_-)=2\c1  B_-\wedge M .\nonumber\\[8pt]
&&\mbox{Non-standard subfamily:}\quad  \c1\nonumber\\[2pt]
&&r=\c1 N\wedge M  \nonumber\\[2pt]
&&\delta(N)=0\qquad \delta(M)=0\label{cb}\\
&&\delta(A_+)=-\c1  A_+\wedge M\qquad 
\delta(A_-)=\c1A_-\wedge M\nonumber\\[2pt]
&&\delta(B_+)=-2\c1  B_+\wedge M\qquad 
\delta(B_-)=2\c1  B_-\wedge M .
\nonumber
\eea

 \noindent
$\bullet$ {\it Type II: $A_+$ primitive}.  If we set
$\delta(A_+)=0$  the initial free parameters are:
$\a1$, $\a3$, $\a4$, $\a5$, $\b3$, $\c1$ and $\c2=-\c1$. The
relations (\ref{bbf})--(\ref{bbh}) reduce to a single
equation 
$\a1\a4+\a5\c1=0$, and the Schouten bracket characterizes  the standard
and non-standard subfamilies by means of the term $\a1\b3-\c1^2$:

\noindent
Standard subfamily: $\a1$, $\a3$, $\a4$, $\a5$, $\b3$, $\c1$ with 
$\a1\a4+\a5\c1=0$,  $\a1\b3-\c1^2\ne 0$.

\noindent
Non-standard subfamily: $\a1$, $\a3$, $\a4$, $\a5$, $\b3$, $\c1$ with 
$\a1\a4+\a5\c1=0$,  $\a1\b3-\c1^2= 0$.

The structure for both subfamilies of bialgebras turns out to be: 
\bea 
&&\!\!\!\!\!\!\!\!\! r=\a1 N\wedge A_+  +
\a3 A_+\wedge M    + \a4 B_+\wedge M +
\a5 A_+\wedge B_+ \cr
&&\!\!\!\!\!\!\!\!\! \qquad     +
\b3 A_-\wedge M  +\c1 (N\wedge M - A_+\wedge A_-)
\nonumber\\[2pt]
&&\!\!\!\!\!\!\!\!\! \delta(\aa)=\a1 \aa\wedge\ap  + \a3
\ap\wedge
\bb + 2 \a4
\bp\wedge\bb + 3 \a5 \ap \wedge \bp  -  \b3 \am\wedge \bb  
\nonumber\\[2pt] 
&&\!\!\!\!\!\!\!\!\! \delta(\ap)=   0\qquad
\delta(\bb)= 0
\label{cc}\\ 
&&\!\!\!\!\!\!\!\!\! \delta(\am)=\a1 (\aa\wedge \bb -
\ap\wedge
\am )    + 2 \a4 \ap\wedge\bb - \a5
\bp\wedge\bb   + 2\c1 \am\wedge \bb   \nonumber\\[2pt]
&&\!\!\!\!\!\!\!\!\!\delta(\bp)=   2 \a1 \ap\wedge \bp   -
2\b3\ap
\wedge \bb  -2\c1\bp\wedge\bb\nonumber\\[2pt]
&&\!\!\!\!\!\!\!\!\!\delta(\bm)=     2 \a1 (\aa\wedge\am -
\ap\wedge \bm  )  
+ 2\a3\am \wedge \bb +4\a4\aa\wedge \bb \cr
&&\!\!\!\!\!\!\!\!\!\qquad\qquad - 2 \a5 (2\aa \wedge \ap
- \am \wedge \bp - \ap \wedge \bb)
+2\c1\bm\wedge\bb .
\nonumber
\eea

\noindent
$\bullet$ {\it Type III: $B_+$ primitive}.
Finally, the condition $\delta(B_+)=0$ implies that the  initial free
parameters are:
$\a1$, $\a2$, $\a3$, $\a4$, $\a5$, $\c2$ and $\b6=\a1$.
The equations (\ref{bbf})--(\ref{bbh}) lead to   $\a1=0$ and 
$\a2\a3-\a5\c2=0$. Hence from (\ref{bj}) we find that standard solutions 
correspond to considering the set of parameters
$\{\a2,\a3,\a4,\a5=\frac{\a2\a3}{\c2},\c2\ne 0\}$:
\bea 
&&\mbox{Standard subfamily:}\quad  \a2,\ 
 \a3,\ \a4,\ \c2\ne 0 \nonumber\\[2pt]
&&r=\a2  N\wedge B_+  + \a3  A_+\wedge M +\a4  B_+\wedge M
 +\frac{\a2\a3}{\c2}  A_+\wedge B_+ + \c2
A_+\wedge A_-   \nonumber\\[2pt]
&&\delta(N)=2 \a2 N \wedge B_+ +\a3  A_+ \wedge M
+ 2\a4  B_+\wedge M +3\frac{\a2\a3}{\c2}  A_+\wedge B_+\cr
&&\delta(A_+)=\a2  B_+\wedge A_+ - \c2  A_+
\wedge M\label{ce} \\
&&\delta(A_-)=-
\a2  (B_+\wedge A_-  + 2  A_+\wedge N) +2 \a4  A_+  \wedge M
  -  \frac{\a2\a3}{\c2}  B_+\wedge M\cr
&&\qquad\qquad  -\c2  A_- \wedge M\nonumber\\[2pt]
&&\delta(B_+)=0\qquad \delta(M)=0\nonumber\\[2pt]
&&\delta(B_-)=-2\a2  B_+ \wedge B_-
+ 2\a3  A_- \wedge M +2(\a2 +2\a4)
N\wedge M\cr
&&\qquad\qquad - 2\frac{\a2\a3}{\c2}(2 N\wedge A_+ -  A_- \wedge B_+ - 
A_+\wedge M).\nonumber
\eea
The non-standard subfamily corresponds to  taking $\{\a2,\a3,\a4,\a5\}$
together with the relation  $\a2\a3=0$ which implies that either
$\a2$ or $\a3$ is equal to zero. However if  we set $\a2=0$ then
$\delta(A_+)=0$ and we are within the above non-standard type II;
therefore we discard it and only consider the case $\a3=0$.
\bea 
&&\mbox{Non-standard subfamily:}\quad  \a2,\ 
\a4,\  \a5\nonumber\\[2pt]
&&r=\a2  N\wedge B_+ + \a4 B_+\wedge M 
+\a5  A_+\wedge B_+   \nonumber\\[2pt]
&&\delta(N)=2 \a2 N \wedge B_+ 
+ 2\a4  B_+\wedge M +3\a5  A_+\wedge B_+\label{cf}\\
&&\delta(A_+)=\a2  B_+\wedge A_+\nonumber\\[2pt]
&&\delta(A_-)=-\a2  (B_+\wedge A_-  + 2  A_+\wedge N) +2 \a4  A_+ 
\wedge M
  -  \a5  B_+\wedge M\nonumber\\[2pt]
&&\delta(B_+)=0\qquad   \delta(M)=0\nonumber\\[2pt]
&&\delta(B_-)=-2\a2  B_+ \wedge B_-
+2(\a2 +2\a4)N\wedge M\cr
&&\qquad\qquad
- 2\a5(2 N\wedge A_+ - A_- \wedge B_+ -  A_+\wedge
M).
\nonumber
\eea

In what follows we will study the quantum deformations of  two specific 
bialgebras of non-standard type  with either $A_+$ or $B_+$ as primitive
generators. The former contains a quantum harmonic  oscillator $h_4$
subalgebra  while the latter includes a quantum $gl(2)$ subalgebra.


\section{The quantum two-photon algebra $U_{\a1}(h_6)$}

We consider the  bialgebra belonging to the non-standard  subfamily of
type II with $\a3=\a4=\a5=\b3=\c1=0$ and $\a1$ as the only free
parameter. Thus this one-parameter two-photon bialgebra can be written as
\be
\begin{array}{l}
r=\a1 N\wedge A_+  \cr
\delta(\ap)=0\qquad \delta(\bb)=0\cr
\delta(\aa)=\a1\aa\wedge\ap\qquad 
\delta(\bp)=-2\a1\bp\wedge\ap\cr
\delta(\am)=\a1(\am\wedge\ap + \aa\wedge\bb)\cr
\delta(\bm)=2\a1(\bm\wedge\ap + \aa\wedge\am). 
\end{array}
\label{db}
\ee
In order to construct its corresponding quantum  deformation 
$U_{\a1}(h_6)$, it is necessary to obtain a  homomorphism, the {\em
coproduct},  $\Delta: U_{\a1}(h_6)\to U_{\a1}(h_6)\otimes U_{\a1}(h_6)$,
verifying the coassociativity condition
$(\Delta\otimes\mbox{id})\Delta=\Delta(\Delta\otimes\mbox{id})$.  It
turns out to be
\cite{tufotona}:
\be
\begin{array}{l}
\Delta(\ap)=1\otimes \ap + \ap \otimes 1\qquad
\Delta(\bb)=1\otimes \bb + \bb \otimes 1\cr
\Delta(\aa)=1\otimes \aa + \aa \otimes e^{\a1\ap}
\qquad
\Delta(\bp)=1\otimes \bp + \bp \otimes e^{-2\a1\ap}\cr
\Delta(\am)=1\otimes \am + \am \otimes e^{\a1\ap}
+ \a1\aa \otimes e^{\a1\ap} \bb\cr
\Delta(\bm)=1\otimes \bm + \bm \otimes e^{2\a1\ap}
-\a1\am\otimes e^{\a1\ap} \aa \cr
\qquad\qquad\quad 
+ \a1\aa \otimes e^{\a1\ap}(\am-\a1 \bb\aa) . 
\end{array}
\label{dc}
\ee
The compatible deformed commutation rules are obtained by
imposing $\Delta$ to be a homomorphism of $U_{\a1}(h_6)$: $\Delta
([X,Y])=[\Delta(X),\Delta(Y)]$; they are
\bea
&& [\aa,\ap]=\frac {e^{\a1\ap}-1}{\a1} \qquad   [\aa,\am]=-\am \qquad
 [\am,\ap]=\bb  e^{\a1\ap} \cr
&& [\aa,\bp]=2\bp \qquad   [\aa,\bm]=-2\bm - \a1 \am \aa
\qquad [\bb,\,\cdot\,]=0\nonumber\\[2pt]
&&[\bm,\bp]= 2(1+  e^{-\a1\ap})\aa + 2\bb - 2 \a1 \am\bp  \label{dd}\\
&& [\ap,\bm]=-(1+  e^{\a1\ap}) \am + \a1 e^{\a1\ap}\bb\aa\qquad 
[\ap,\bp]=0\cr
&&  [\am,\bp]=2\frac {1-e^{-\a1\ap}}{\a1} 
  \qquad   [\am,\bm]=-\a1\am^2 . 
\nonumber 
\eea
The associated universal quantum R-matrix, which satisfies the quatum
YBE, reads
\be
{\cal R}=\exp\{-\a1\ap\otimes\aa\}\exp\{\a1\aa\otimes\ap\}  .
\label{de}
\ee
Note that the underlying cocommutator  is related to the first  order in
$\a1$ of the coproduct  by means of $\delta=(\Delta - \sigma\Delta)$
where
$\sigma(X\otimes Y)=Y\otimes X$; the limit $\a1\to
0$ of (\ref{dd}) leads to the $h_6$ Lie brackets (\ref{aa});  and the
first order in $\a1$ of ${\cal R}$ corresponds to the classical
$r$-matrix (\ref{db}). We remark that the generators
$\{\aa,\ap,\am,\bb\}$ give rise to a non-standard quantum harmonic
oscillator algebra  $U_{\a1}(h_4)\subset U_{\a1}(h_6)$ \cite{osc}.

On the other hand, a deformed one-boson realization of  $U_{\a1}(h_6)$
is given by
\bea
&&\aa= \frac{e^{\a1\aap}-1}{\a1}\,\aam \qquad
  \ap=\aap  \qquad 
  \am=  e^{\a1\aap}\aam \cr
&& \bp=\left(\frac{1- e^{-\a1\aap} }{\a1}\right)^2  \qquad 
\bm=e^{\a1\aap}\aam^2  \qquad  \bb= 1  .
\label{dg}
\eea
Notice that the classical identifications  $B_+=A_+^2$ and $B_-=A_-^2$
are no longer valid in the quantum case.  The action of the generators of
$U_{\a1}(h_6)$ on the number states
$\{|m\rangle\}_{m=0}^{\infty}$  is
\bea
&&\ap|m\rangle =\sqrt{m+1}\, |m+1\rangle  \qquad
\bb|m\rangle =  |m\rangle   \cr
&&\am|m\rangle = \sqrt{m}\,
|m-1\rangle + m \sum_{k=0}^\infty
 \frac{{\a1}^{k+1} }{(k+1)!} \, \sqrt{\frac {(m+k)!}{m!} }
 \, |m+k\rangle  \cr
&&\aa|m\rangle =m\, |m\rangle + m \sum_{k=1}^\infty  \frac{{\a1}^k
}{(k+1)!} \, \sqrt{\frac {(m+k)!}{m!} }
 \, |m+k\rangle  \cr
&&\bp|m\rangle =\sqrt{(m+1)(m+2)}\,|m+2\rangle \cr
&&\qquad \qquad + \sum_{k=1}^\infty
(-2+2^{k+2}) \frac{(-{\a1})^k }{(k+2)!} \, \sqrt{\frac {(m+k+2)!}{m!} }
 \, |m+k+2\rangle  \cr
&&\bm|m\rangle =\sqrt{m (m-1)}\, |m-2\rangle +
{\a1} \sqrt{m} (m-1)\, |m-1\rangle \cr
&&\qquad \qquad + m (m-1) \sum_{k=0}^\infty  \frac{{\a1}^{k+2}
}{(k+2)!} \, \sqrt{\frac {(m+k)!}{m!} }
 \, |m+k\rangle  .\label{dh}
\eea

The deformed boson realization (\ref{dg}) can be 
translated into differential operators acting on the the space of
entire analytic functions $f(\alpha)$, that is, a deformed Fock--Bargmann
representation which is given by
\bea
&& \aa= \frac{e^{{\a1}\alpha}-1}{{\a1}}\,\frac{d}{d\alpha}  \qquad
  \ap=\alpha   \qquad 
  \am=  e^{\a1\alpha}\,\frac{d}{d\alpha}    \cr
&& \bp=\left(\frac{1- e^{-\a1\alpha }}{\a1}\right)^2   \qquad 
\bm=e^{\a1\alpha}\,\frac{d^2}{d\alpha^2} \qquad  \bb= 1 .
\label{di}
\eea
Hence the relation (\ref{ah}) provides the following   differential
equation that characterizes the   deformed two-photon algebra
eigenstates :
\bea
&&\beta_2 e^{\a1\alpha} \frac{d^2f}{d\alpha^2}+\left(\beta_1
\frac{e^{{\a1}\alpha}-1}{{\a1}}  +\beta_4  e^{\a1\alpha}\right) \frac{d
f}{d\alpha} \cr
 &&\qquad\quad +\left(\beta_3\left(\frac{1- e^{-\a1\alpha
}}{\a1}\right)^2 +\beta_5\alpha -\lambda\right)f=0 .
\label{dj}
\eea
The particular equation with $\beta_2=\beta_3=0$  is associated to the
quantum oscillator subalgebra $U_{\a1}(h_4)$ and it would give deformed  
one-photon coherent states; the case
$\beta_4=\beta_5=0$ corresponds to the $gl(2)$ sector  which
is not a quantum subalgebra (see the coproduct (\ref{dc})). 

We stress the relevance of the coproduct in  order to construct tensor
product representations of the two-photon generators (\ref{di}).
Finally, we remark that the limit $\a1\to 0$ of all the
above expressions gives rise to their classical version
presented in the Introduction.


\section{The quantum two-photon algebra  $U_{\a2}(h_6)$}

Let us consider now the non-standard bialgebra of  type III with
$\a4=\a5=0$ and
$\a2$ as a free parameter:
\be
\begin{array}{l}
r= \a2 \aa\wedge \bp\cr
 \delta(\bp)=0\qquad \delta(\bb)=0\cr
 \delta(\aa)=2\a2\, \aa\wedge \bp\qquad 
 \delta(\ap)=-\a2\, \ap\wedge \bp\cr
 \delta(\am)=\a2\, (\am\wedge \bp + 2 \aa\wedge \ap)\cr
 \delta(\bm)=2\a2\, (\bm\wedge \bp + \aa\wedge \bb) . 
\end{array}
\label{ea}
\ee
The resulting  coproduct, commutation rules and universal  $R$-matrix of
the quantum algebra  $U_{\a2}(h_6)$ read \cite{tufotonb}:
\be
\begin{array}{l}
\Delta(\bp)=1\otimes \bp + \bp \otimes 1 \qquad
\Delta(\bb)=1\otimes \bb + \bb \otimes 1\cr
\Delta(\aa)=1\otimes \aa + \aa \otimes e^{2\a2\bp}
\qquad \Delta(\ap)=1\otimes \ap + \ap \otimes
e^{-\a2\bp} \cr
\Delta(\am)=1\otimes \am + \am \otimes e^{\a2\bp}
+ 2\a2\aa \otimes e^{2\a2\bp} \ap\cr
\Delta(\bm)=1\otimes \bm + \bm \otimes e^{2\a2\bp}
+ 2\a2\aa \otimes e^{2\a2\bp} \bb\cr
\end{array}
\label{eb}
\ee
\bea
&& [\aa,\ap]=\ap \qquad   [\aa,\am]=-\am \qquad
 [\am,\ap]=\bb   \cr
&& [\aa,\bp]=\frac {e^{2\a2\bp}-1}{\a2} \qquad   [\aa,\bm]=-2\bm - 4\a2
{\aa}^2 \cr
&&[\bm,\bp]= 4\aa + 2\bb e^{2\a2\bp} \qquad  [\bb,\,\cdot\,]=0 
\label{ec}\\
 && [\ap,\bm]=-2\am + 2\a2(\aa \ap + \ap \aa) \qquad
[\ap,\bp]=0\cr
&&  [\am,\bp]=2 e^{2\a2\bp} \ap
  \qquad   [\am,\bm]=-2\a2(\aa \am + \am \aa) 
\nonumber
\eea
 \be
{\cal R}=\exp\{-\a2\bp\otimes\aa\}\exp\{\a2\aa\otimes\bp\} .
\label{ed}
\ee
Notice that the generators $\{N,B_+,B_-,M\}$  close a non-standard
quantum
$gl(2)$ algebra \cite{gl} such that  $U_{\a2}(gl(2))\subset  
U_{\a2}(h_6)$, while the  oscillator  algebra $h_4$ is 
preserved as an undeformed subalgebra  only at the level of
commutation relations.

A one-boson representation of $U_{\a2}(h_6)$   is given by: 
\bea
&&\bp=\aap^2\qquad \bb=1\qquad
 \aa=\frac{e^{2\a2\aap^2}-1}{2\a2\aap}\,\aam \cr
&& \ap=\left(\frac{1-e^{-2\a2\aap^2}}{2\a2}\right)^{1/2}\qquad 
 \am=\frac{e^{2\a2\aap^2}}{\aap}
\left(\frac{1-e^{-2\a2\aap^2}}{2\a2}\right)^{1/2} \!\! \aam \cr
&& \bm =\left( \frac{e^{2 \a2 \aap^2}-1}{2 \a2 \aap^2 }   \right)   
\aam^2 + 
\left(\frac {e^{2 \a2 \aap^2}}{\aap}+\frac {1-e^{2 \a2 \aap^2}} {2 \a2
\aap^3}\right)  
\aam .
\label{ee}
\eea
 We remark that, although
(\ref{ec}) presents a non-deformed oscillator subalgebra,  the
representation (\ref{ee}) includes strong deformations in terms of the
boson operators. The corresponding Fock--Bargmann representation 
adopts the following  form:
\bea
&&\bp=\alpha^2\qquad \bb=1\qquad
 \aa=\frac{e^{2\a2\alpha^2}-1}{2\a2\alpha}\,\frac{d}{d\alpha} \cr
&& \ap=\left(\frac{1-e^{-2\a2\alpha^2}}{2\a2}\right)^{1/2}\qquad 
 \am=\frac{e^{2\a2\alpha^2}}{\alpha}
\left(\frac{1-e^{-2\a2\alpha^2}}{2\a2}\right)^{1/2} \!\! 
\frac{d}{d\alpha} \cr && \bm =\left( \frac{e^{2 \a2 \alpha^2}-1}{2 \a2
\alpha^2 }   \right)   
\frac{d^2}{d\alpha^2}  + 
\left(\frac {e^{2 \a2 \alpha^2}}{\alpha}+\frac  {1-e^{2 \a2 \alpha^2}}
{2 \a2
\alpha^3}\right)  
\frac{d}{d\alpha} .
\label{ef}
\eea
Therefore if  we substitute these operators in the equation of the
two-photon algebra eigenstates    
(\ref{ah}) we obtain the  differential equation:
\bea
&&
\!\!\!\!\!\!\!\!\!\!\!\!
\beta_2 \left( \frac{e^{2 \a2 \alpha^2}-1}{2 \a2 \alpha^2 }  
\right)  
\!\frac{d^2f}{d\alpha^2}
+\left(\beta_1 \frac{e^{2\a2\alpha^2}-1}{2\a2\alpha}
 + \beta_4
\frac{e^{2\a2\alpha^2}}{\alpha}
\left(\frac{1-e^{-2\a2\alpha^2}}{2\a2}\right)^{1/2} \right)\!
\frac{d f}{d\alpha} \cr
&&
\!\!\!\!\!\!\!\!\!\!\!\!
+\beta_2  \left(\frac {e^{2 \a2 \alpha^2}}{\alpha}+\frac
{1-e^{2 \a2
\alpha^2}} {2 \a2
\alpha^3}\right)\! \frac{d
f}{d\alpha} +\left(\beta_3 \alpha^2    
+\beta_5\left(\frac{1-e^{-2\a2\alpha^2}}{2\a2}\right)^{1/2}
\!\! -\lambda\right)\! f=0 .\cr
&&
\!\!\!\!\!\!\!\!\!\!\!\!
\label{eg}
\eea
If we set $\beta_4=\beta_5=0$,  then we obtain an equation associated to
the quantum subalgebra $U_{\a2}(gl(2))$, while the case
$\beta_2=\beta_3=0$ corresponds to the harmonic oscillator sector.
Note that in the limit $\a2\to 0$ we recover the classical
two-photon structure.



To end with, it is remarkable that we can make use of the two-photon
bialgebra automorphism  defined by (\ref{bg})
and (\ref{bh}) in order to obtain from $U_{\a1}(h_6)$ and 
$U_{\a2}(h_6)$ two  other (algebraically equivalent) quantum
deformations of $h_6$, namely $U_{\b1}(h_6)$ and 
$U_{\b2}(h_6)$, but now with $A_-$ and
$B_-$ as the primitive generators, respectively. However at a
representation level, the physical implications are rather
different. If, for instance, $A_-$ is a primitive generator
instead of $A_+$ we would obtain a deformed 
Fock--Bargmann representation with terms as
$\exp(\b1\frac{d}{d\alpha})$ (instead of $e^{\a1\alpha}$), 
giving rise to a differential-difference realization
\cite{tufotona}. Therefore, quantum $h_6$ algebras with
either $\ap$ or $\bp$ primitive  would
originate a class of smooth deformed states, while
those with either $\am$ or $\bm$ primitive  will be
linked to a set of states including some intrinsic
discretization.


{\section*{Acknowledgments}}
\noindent
 A.B. and
F.J.H. have been partially supported by DGICYT (Project  PB94--1115)
from the Ministerio de Educaci\'on y Cultura  de Espa\~na and by Junta
de Castilla y Le\'on (Pro\-jects CO1/396 and CO2/297). P.P. has been
supported by a fellowship from AECI, Spain.


\end{document}